\documentclass[11pt,twoside]{amsart}
\usepackage{amsmath}
\usepackage{amsthm}
\usepackage{amsfonts}
\usepackage{amssymb}
\usepackage[all]{xy}
\usepackage{alltt}

\theoremstyle{plain}
\newtheorem{prop}{Proposition}[section]
\newtheorem{conjecture}[prop]{Conjecture}
\newtheorem*{thmA}{Theorem A}

\newtheorem{thm}[prop]{Theorem}
\newtheorem{cor}[prop]{Corollary}
\theoremstyle{definition}
\newtheorem{defn}[prop]{Definition}

\newtheorem{remark}[prop]{Remark}
\newtheorem{example}[prop]{Example}

\newcommand{\ints}{\ensuremath{\mathbb{Z}}}
\newcommand{\intmod}[1]{\ensuremath{\mathbb{Z}/{(#1)}}}

\newcommand{\rats}{\ensuremath{\mathbb{Q}}}

\newcommand{\proj}[1]{\ensuremath{\mathbb{P}^{#1}}}

\newcommand{\n}{\noindent}

\newcommand{\lspace}{\hspace{3pt}}

\DeclareMathOperator{\Hom}{Hom}

\DeclareMathOperator{\Ext}{Ext}

\DeclareMathOperator{\Tor}{Tor}

\newcommand{\smsh}{\ensuremath{{\scriptstyle \wedge}}}             

\newcommand{\wdg}{\ensuremath{{\scriptstyle \vee}}}

\newcommand{\iprod}[1]{\ensuremath{\underset{\scriptscriptstyle #1}{\Pi}}}

\DeclareMathOperator{\im}{im}
\DeclareMathOperator{\coker}{coker}
\DeclareMathOperator{\coim}{coim}

\newcommand{\map}[3]{\ensuremath{#1 : #2 \longrightarrow #3}}

\newcommand{\stabcat}{\ensuremath{\mathcal{S}}}

\DeclareMathOperator{\ann}{ann}
\DeclareMathOperator{\Rann}{Rann}
\DeclareMathOperator{\Lann}{Lann}

\newcommand{\gh}{\ensuremath{\mathcal{GH}}}

\newcommand{\thicksub}[1]{\ensuremath{{{\bf thick}\langle#1\rangle}}}

\newcommand{\dr}{\ensuremath{{\mathcal{D}(R)}}}

\DeclareMathOperator{\spec}{Spec}

\newcommand{\qcoh}[1]{\ensuremath{\mathfrak{Qcoh}(#1)}}
\newcommand{\dqcoh}[1]{\ensuremath{\mathcal{D}(\qcoh{#1})}}

\newcommand{\modcat}[1]{\ensuremath{#1-{\bf Mod}}}

\newcommand{\pistar}{\ensuremath{\pi_*\hspace{1pt}}}

\newcommand{\chain}[1]{\ensuremath{{#1}_\text{\textbullet}}}

\begin{document}  

\title{The Generating Hypothesis in the Derived Category of $R$-modules}
\date{\today}

\author{Keir H. Lockridge} 
\address {Department of Mathematics \\
          University of Washington \\
          Seattle, WA  98195 }
\email{lockridg@math.washington.edu}

\subjclass[2000]{Primary: 55p42; Secondary: 18e30}

\begin{abstract} In this paper, we prove a version of Freyd's 
generating hypothesis for triangulated categories:  if 
$\mathcal{D}$ is a cocomplete triangulated category and $S \in 
\mathcal{D}$ is an object whose endomorphism ring is 
graded commutative and concentrated 
in degree zero, then $S$ generates (in the sense of Freyd) the 
thick subcategory determined by $S$ if and only if the 
endomorphism ring of $S$ is von Neumann regular.  As a corollary, 
we obtain that the generating hypothesis is true in the derived 
category of a commutative ring $R$ if and only if $R$ is von 
Neumann regular.  We also investigate alternative formulations of 
the generating hypothesis in the derived category.  Finally, we 
give a characterization of the Noetherian stable homotopy 
categories in which the generating hypothesis is true. 
\end{abstract}

\maketitle

\tableofcontents

\section{Introduction}

In his 1966 paper {\em Stable Homotopy}, Freyd introduces the following conjecture:

\begin{conjecture}[The Generating Hypothesis] If $\map{f}{X}{Y}$ is a morphism of finite spectra and $\pi_* f = 0$, then f is trivial.\label{gh}
\end{conjecture}

\n The conjecture remains open, but progress has been made using the methodology surrounding the Nilpotence Theorem of Devinatz, Hopkins, and Smith (\cite{devinatz-1988}, \cite{hopkins-smith-1998}).  Let $L_1(-)$ denote $E(1)$-localization.  In \cite{devinatz-1990}, Devinatz proves the $E(1)$-approximate generating hypothesis when the target spectrum is the sphere:  if $\pi_* f = 0$, then $L_1 f \simeq 0$.  An axiomatic approach to stable homotopy theory has led to the study of other triangulated categories from a homotopy theoretic point of view (\cite{hovey-axiomatic-1997}), and a natural extension of this study is to try to formulate and prove the generating hypothesis in these structurally similar settings.  Examples of general stable homotopy categories include localizations of the usual stable category, the derived category of a ring, the stable module category (arising in representation theory), the homotopy category of complexes of injective comodules over a Hopf algebra (e.g., the Steenrod algebra), among others.

It is not completely clear how the generating hypothesis should be stated in general stable categories;  one can leave the statement as it is and get something that makes perfect sense, but there are several ways of characterizing the finite spectra from an axiomatic point of view that do not always coincide in general.  Difficulties arise, for example, when the sphere is not a weak generator; i.e., $\pi_*X = 0$ does not necessarily imply that $X$ is trivial.  If such a finite $X$ exists, then \gh \lspace would clearly be false (an example will be given later: the derived category of quasi-coherent sheaves over $\proj{1}$).  One might insist that $[X,f] = 0$ for each weak generator $X$.  

We now give two relevant definitions and formulate the version of
\gh\lspace with which we will work.  Let $\mathcal{D}$ be a
triangulated category, and let $S \in \mathcal{D}$ be a distinguished
object.  The {\em thick subcategory generated by $S$}, $\thicksub{S}$,
is the smallest class of objects in $\mathcal{D}$ that contains $S$
and is closed under suspension, retraction, and cofiber sequences.  In
the usual stable category, this is exactly the finite spectra; in
general, it may not describe the class of small objects.  An object
$Y$ is {\em small} if, whenever the coproduct $\coprod_\alpha
X_\alpha$ exists, the natural map $$\xymatrix{\underset{\alpha}{\bigoplus}
  [Y,X_\alpha] \ar@<1ex>[r] & [Y,\underset{\alpha}{\coprod} X_\alpha]}$$ is
an isomorphism.  Note that 'small' and 'finite' have the same meaning
in the usual stable category.  

\begin{defn} Let $\mathcal{D}$ be a triangulated category and let $S \in \mathcal{D}$ be a distinguished object.  Write $\pi_*(-)$ for the functor $\Hom^{*}_\mathcal{D}(S,-)$.  The {\em generating hypothesis} is the statement:  If $\map{f}{X}{Y}$ is a morphism of objects in $\thicksub{S}$ and $\pi_* f = 0$, then $f \simeq 0$.\label{maindef}\end{defn}

In a stable homotopy category, $S$ will be the sphere object.  In the derived category of a ring, for example, the sphere is the chain complex with $R$ in degree zero and zero elsewhere.  The following theorem is our main result.

\begin{thmA} Let $\mathcal{D}$ be a triangulated category with arbitrary coproducts, and let $S$ be an object in $\mathcal{D}$ such that $R = \pistar S$ is commutative and concentrated in degree zero. \gh \lspace is true in $\mathcal{D}$ if and only if $R$ is von Neumann Regular. \end{thmA}

\n A ring $R$ is {\em von Neumann regular} if, for every element $x \in R$, there exists an element $y \in R$ such that $xyx = x$.  This class of rings, in the context of noncommutative ring theory, was originally introduced by von Neumann to study operator algebras in functional analysis.  Theorem A has the following corollary.

\begin{cor} Let $R$ be a commutative ring.  \gh \lspace is true in the derived category of $R$-modules if and only if $R$ is von Neumann regular. \end{cor}

The main result actually applies to any triangulated category where $\pistar S$ is graded commutative and concentrated in even degrees. If $\pistar S$ is graded commutative and {\em not} concentrated in degree zero, we prove that if $\pi_* S$ is nonzero in only finitely many degrees or if $\pi_* S$ connective and concentrated in even degrees, then the generating hypothesis is necessarily false.  We also give an application of Theorem A to the derived category of quasicoherent sheaves over certain schemes in $\S \ref{vnrsection}$.

There is a more abstract way to view the content of Theorem A.  An analysis of the proof shows that, for categories covered by the main result, the generating hypothesis is true if and only if the category is trivial in the following sense:  $\thicksub{S}$ must be exactly the collection of retracts of finite wedges of suspensions of $S$.  In ordinary stable homotopy, this condition would imply that every finite spectrum is a wedge of suspensions of the sphere (the main theorem, of course, does not apply in this case).  We conjecture that this characterization is valid for the derived category of any Grothendieck abelian category.  Note that \stabcat\lspace is not of this form - for any object $X$ in the derived category of an abelian category, the cofiber $Y$ of $\map{2}{X}{X}$ has the property that $\map{2}{Y}{Y}$ is trivial.  However, $\map{2}{M(2)}{M(2)}$ is nontrivial in \stabcat, where $M(2)$ is the mod 2 Moore spectrum (I learned this from Neil Strickland).

Unless otherwise indicated, for the remainder of the paper we consider a triangulated category $\mathcal{D} $ and an object $S \in \mathcal{D}$ which has the property that its endomorphism ring $R = \pi_*S = [S,S]_*$ is commutative and concentrated in degree zero.  We further assume that $\mathcal{D}$ has arbitrary coproducts; this guarantees that idempotents split in $\mathcal{D}$ (\cite[3.2]{bokstedt-neeman}).

This paper is organized as follows.  In \S 2, we develop two important necessary conditions for \gh\lspace to hold: the Nilpotence Criterion (\ref{nilcrit}) and the Annihilator Criterion (\ref{anncrit}).  In \S 3, we introduce (von Neumann) regular rings and show that the two criteria from \S 2 are equivalent to regularity.  We then show that regularity is also sufficient.  In \S 4, we explore variants of \gh\lspace by changing its domain of definition.  In $\dr$, for example, we study \gh\lspace when stated for maps of objects in $\thicksub{S/I}$, where $I$ is a finitely generated ideal of $R$.  In \S 5, we show that \gh\lspace holds in a Noetherian stable homotopy category if and only if $\thicksub{S}$ is exactly the collection of retracts of finite wedges of suspensions of $S$.

I would like to thank my advisor, Ethan Devinatz, for many helpful conversations regarding the content and preparation of this paper.  Further, I wish to acknowledge the contribution of the referee, who suggested the given proof of Proposition \ref{vnrsummand}.  The original proof did not apply to noncommutative regular rings, and in its present form, the proposition allows us to make Remark \ref{noncommincursion}.  The referee also suggested the inclusion of Proposition \ref{dcatobs} which is used to establish the noncommutative aspect of Corollary \ref{noethcase}.

\section{Generalities and Criteria}

By way of motivation, let us consider the derived category of a commutative ring \dr, where $S$ is the chain complex with $R$ concentrated in degree zero.  \dr\lspace is a monogenic stable homotopy category (\cite[9.3.1]{hovey-axiomatic-1997}).  If $R$ is a field, then every object in \dr\lspace is equivalent to a wedge of suspensions of the sphere; therefore, \gh\lspace is trivially true.  The following proposition gives an apposite connection between \gh\lspace and direct products.

\begin{prop}  Suppose $R \cong R_1 \times R_2$.  \gh\lspace is true in $\dr$ if and only if it is true in $\mathcal{D}(R_1)$ and $\mathcal{D}(R_2)$. \label{sums}\end{prop}
\begin{proof}
This is essentially a consequence of the fact that every $R$-module is the direct sum of an $R_1$-module and an $R_2$-module and every $R$-module map decomposes similarly.
\end{proof}

\begin{remark}Hence, \gh\lspace is true in \dr\lspace whenever $R$ is a finite product of fields.  In particular, \gh\lspace is true for $\intmod{n}$, provided $n$ is square free.  This condition is also necessary.  In \cite{me} we prove that if \gh\lspace is true in a category where $S$ is connective and $\pi_0 S$ is projective-free (meaning that every projective $\pi_0 S$ module is free), then $\pi_* S$ is either a field or totally non-coherent.  A graded ring is totally non-coherent if no proper, nonzero, finitely generated ideal is finitely presented.  Hence, \gh\lspace is false for $\intmod{p^n}$ when $n > 1$ since $\intmod{p^n}$ is local and therefore projective-free but obviously neither a field nor totally non-coherent.  We also give a specific counterexample in Example \ref{zmodn}.  \label{product_fields}\end{remark}

The next proposition is based upon Devinatz' approach to \gh\lspace in \cite{devinatz-1990}, where he proves the $E(1)$-approximate version when the target is the sphere.  The Spanier-Whitehead duality pairing $\xymatrix@1{DX\smsh X \ar[r] & S}$ induces a map
\begin{equation}
\xymatrix{\pi_{-i}(DX)\otimes_R \pi_i(X) \ar[r] & R.} \label{pairing}
\end{equation}
Let $M$ be an $R$-module.  After tensoring with $M$, we may examine the adjoint map
\begin{equation}
\xymatrix{\pi_{-i}(DX) \otimes_R M \ar[r] & \Hom_R(\pi_i(X), M).} \label{trans}
\end{equation}
If $M$ is both flat and injective, then this is a natural transformation of cohomology theories.  Since it is an isomorphism when $X = S$, it is an isomorphism for all $X \in \thicksub{S}$.  Now suppose there is an injection $\xymatrix@1{R \ar[r] & M}$ and $M/R$ is flat.  Then 
\[ \xymatrix{\pi_{-i}(DX) \ar[r] & \pi_{-i}(DX) \otimes_R M} \] 
is also injective since the failure of this map to be injective is measured by \[\Tor^1_R(\pi_{-i}(DX), M/R),\] which is trivial.  Now consider a degree zero map $\map{f}{X}{S}$ where $X \in \thicksub{S}$.  The map $f$ corresponds to a map $\map{f'}{S}{DX}$.  If $f$ is nontrivial, then so is $f'$; $f'$ is also nontrivial in $\pi_0(DX) \otimes_R M$.  Since (\ref{trans}) is an isomorphism and the pairing (\ref{pairing}) corresponds to composition, there is a map $\map{g}{S}{X}$ such that $fg$ is nontrivial.  Hence, $\pi_* f \neq 0$.  We have proved:

\begin{prop}  If there is an $R$-module $M$ that is flat and injective and if there is an injective map $\xymatrix@1{R \ar[r] & M}$ such that $M/R$ is flat, then \gh\lspace for the target spectrum $S$ is true in $\dr$. \end{prop}

\n For example, if $R$ is a von Neumann regular ring (see \S 3), then every $R$-module is flat (\cite[4.2.9]{weibel-introduction-1994}); taking $M$ to be the injective hull of $R$, we obtain that \gh\lspace for the target $S$ is true in \dr.  Also, if $R$ is self injective, then \gh\lspace for the target $S$ is true in $\dr$.  As we will see in the next example, however, there are self-injective rings $R$ for which the general form of \gh\lspace is false in $\dr$.

\begin{example} $\intmod{n}$ is self injective for all $n$.  Hence, \gh\space for the target $S$ is true in $\mathcal{D}(\intmod{p^2})$.  We have observed, however, that \gh\lspace is false in $\mathcal{D}(\intmod{p^2})$ in general.  Here is a counterexample:  consider the map of chain complexes
\begin{equation*}
\xymatrix{0 & 0  \\
		\intmod{p^2} \ar[u] \ar[r]^p & \intmod{p^2} \ar[u] \\
		\intmod{p^2} \ar[u]_-p \ar[r]^-0 & \intmod{p^2} \ar[u]^p\\
		0 \ar[u] & 0 \ar[u]  \\
		&&\\
		X \ar[r]^-f & X.}
\end{equation*}
Though $\pi_* f = 0$, $f$ is nontrivial; a null-chain-homotopy of $f$ would correspond to an element $s \in \intmod{p^2}$ such that $p = ps$ and $ps = 0$, implying $p = 0$.\label{zmodn}\end{example}

Our next proposition gives a general reason for the failure of \gh\lspace in the above example.  In this proposition, we return to the general situation, where $\mathcal{D}$ is a triangulated category with distinguished object $S$ such that $\pi_* S$ is commutative and concentrated in degree zero.

\begin{prop}[Nilpotence Criterion] If \gh\lspace is true in $\mathcal{D}$, then $R$ contains no nonzero nilpotent elements.\label{nilcrit} \end{prop}
\begin{proof} For now, we do not assume that $R$ is commutative.  Let $g \in R$ be nonzero nilpotent; replacing $g$ with some power if necessary, we may assume $g^2 = 0$.  Consider the following commutative diagram, where the rows are cofiber sequences.  The map $h$ is any map that makes the diagram commute.
\begin{equation*}
\xymatrix{S \ar[r]^-g \ar[d]^-0& S \ar[r]^\phi \ar[d]^{g} & Y \ar[r]^-\delta \ar@{-->}[d]^-h & \Sigma S \ar[d]^0\\
S \ar[r]^-g & S \ar[r]^-\phi & Y \ar[r]^-\delta & \Sigma S}
\end{equation*}
Notice that $h$ factors through both $S$ and $\Sigma S$; since $\pi_*S$ is concentrated in one degree, we must have $\pi_* h = 0$.  Hence, any $h$ that fills in the above diagram of cofiber sequences must be trivial by \gh.  Now, since $g^2 = 0$, $g = f\phi$.  One may take $h$ to be $\phi f$, and consequently $\phi f \simeq 0$.  Hence, $f = gk$.  Let $x = k\phi$;  then, $g = f \phi = gx$ and  $xg = k\phi g = 0$.  If $g$ lies in the center of $R$, then we have a contradiction.  In particular, if $R$ is commutative, then it contains no nontrivial nilpotent elements.
\end{proof}

\begin{remark} For this remark, assume $\pistar S$ is graded commutative but not necessarily concentrated in degree zero.  If the map $g$ in the proof has degree $k$, then one obtains that $h$ factors through both $\Sigma^{-k} S$ and $\Sigma^{k+1} S$.  Hence, the Nilpotence Criterion holds whenever $\pi_*S$ is concentrated in even degrees.  Notice also that if \gh\lspace is true, then a nonzero element $g$ of degree $k$ cannot be nilpotent if $\pi_*S$ is zero outside a range of $2k+1$ consecutive degrees.  In particular, if $\pi_i S$ is zero for all but finitely many $i$ and is not concentrated in degree zero, then \gh\lspace is false. \label{evendegree}\end{remark}

In \cite{freyd-1966}, Freyd proves that if the generating hypothesis is true in the stable category of spectra, then $\ann \ann (x) = (x)$ for all $x \in \pi_* S$.  We now give a generalization of this result for $\mathcal{D}$.  We have so far assumed that $\pi_*S$ is commutative and concentrated in degree zero; the following proposition is true without these assumptions.  We write $(x)_L$ for the left ideal generated by $x$ and $(x)_R$ for the right ideal generated by $x$.  We also define the right annihilator of a left ideal $I_L$ and the left annihilator of a right ideal $I_R$ by $\Rann I_L = \{x \in R \lspace|\lspace ix = 0 \lspace\text{for all}\lspace i \in I_L\}$, and $\Lann I_R = \{x \in R \lspace|\lspace xi = 0 \lspace\text{for all}\lspace i \in I_R\}$.  $\Rann I_L$ is a right ideal and $\Lann I_R$ is a left ideal.

\begin{prop}[Annihilator Criterion] If \gh\lspace is true in $\mathcal{D}$, then for any $f \in R$,  $\Lann \Rann (f)_L = (f)_L$.\label{anncrit}\end{prop}
\begin{proof}
Let $k$ be the degree of $f$.  The containment $(f)_L \subseteq \Lann \Rann (f)_L$ is always true.  Take any $g \in \Lann \Rann (f)_L$.  Consider the cofiber sequence
\begin{equation*}
\xymatrix{\Sigma^k S \ar[r]^-f & S \ar[r]^-\phi & Y \ar[r]^-\delta &\Sigma^{k+1} S.}
\end{equation*}
Since $\im \pi_* \delta = \ker \pi_* f = \Rann (f)_L$, we have that $\pi_*(g\delta) = 0.$  By \gh\lspace, $g\delta \simeq 0$.  Hence, $g = \psi f$ and so
\begin{equation*}
(f)_L = \Lann \Rann (f)_L
\end{equation*}
as desired.
\end{proof}

\begin{remark} One could just as easily prove that $\Rann \Lann (f)_R = (f)_R$.\end{remark}

\begin{example} Using the Annihilator Criterion, we give a family of examples of totally non-coherent local (and hence projective-free) rings for which \gh\lspace in \dr\lspace is still false.  Let $k$ be a field and consider $R = k[x_1, x_2, \dots]/(x_ix_j, i \neq j, x_i^n, i \geq 1)$ or $k[x_1, x_2, \dots]/(x_ix_j, i \neq j)$.  $R$ is totally non-coherent (exercise).  However,  $(x_1) \not\subseteq (x_1 + x_2)$ though $\ann (x_1 + x_2) \subseteq \ann (x_1)$.\end{example}

The above two necessary conditions, as we shall see, are also sufficient.

\section{Von Neumann Regular Rings \label{vnrsection}}

There is a general term for a ring that satisfies both the Annihilator Criterion and the Nilpotence Criterion.  A ring (not necessarily commutative) $R$ is {\em (von Neumann) regular} if, for every $x \in R$, there exists $y \in R$ such that $xyx = x$.  For commutative rings, this is the same as the requirement that $(x^2) = (x)$ for all $x \in R$, which is in turn equivalent to the condition that every principal ideal be generated by an idempotent: if there exists an element $s \in R$ such that $sx^2 = x$, then $(x) = (sx)$ and $sx$ is idempotent.  \cite{goodearl-neumann-1991} is a general reference for rings of this type.  One can say something concrete about the class of commutative regular rings.  The following two propositions must be well known; we include proofs for the reader's convenience.

\begin{prop} If $R$ is commutative von Neumann regular, then it is a subring of a direct product of fields.\label{vnrchar}\end{prop}
\begin{proof} First, we show that every prime ideal of $R$ is maximal.  Let $P \subseteq R$ be a prime ideal.  For any $x \in R$, there is some $y \in R$ such that $x^2 y = x$ by regularity.  Hence, $x(xy - 1) = 0$.  Since $R/P$ is an integral domain, either $x$ is zero mod $P$ or $xy \equiv 1$ mod $P$.  Therefore $R/P$ is a field and $P$ is a maximal ideal.  Next, we observe that the nilradical $\mathcal{N}(R)$ of $R$ is trivial: $\mathcal{N}(R)$ cannot contain nonzero idempotent elements.  If $e \in \mathcal{N}(R)$ is idempotent, then since $e$ is contained in every prime ideal, $1-e$ is contained in no prime ideal; therefore, $1-e$ is a unit.  But since $e(1-e) = 0$, $e$ must be zero.  Now consider the map $\map{\phi}{R}{\iprod{\mathfrak{m}} R/\mathfrak{m}}$ defined by $\pi_{\mathfrak{m}}\phi(x) = x$ mod $\mathfrak{m}$.  This is an injective ring homomorphism.
\end{proof}

Any direct product of fields is regular, but the converse of the above Proposition is easily seen to be false; for example, $\ints \subseteq \rats$ is not a regular ring (and, of course, \gh\lspace is false in $\mathcal{D}(\ints)$).

\begin{prop} A commutative ring $R$ is regular if and only if it satisfies both the Annihilator Criterion and the Nilpotence Criterion. \label{regequiv}\end{prop}
\begin{proof}  Suppose $R$ is commutative regular.  In general, if $e$ is idempotent, then $\ann (e) = (1-e)$.  Since each principal ideal $(x)$ is generated by an idempotent, we obtain \[\ann \ann (x) = (x)\] for all $x \in R$.  Further, since $(x) = (x^2)$, it is clear that there are no nonzero nilpotent elements of $R$.

Conversely, suppose the two criteria hold.  Every element of $\ann (x) \cap (x)$ has the property that its square is zero; hence, by the Nilpotence criterion, this intersection is trivial.  If $x^2 y = 0$, then $xy \in \ann (x)$, forcing $xy = 0$.  Thus $\ann (x^2) \subseteq \ann (x)$, and by the Annihilator Criterion, $(x) \subseteq (x^2)$.   Hence $R$ is regular.\end{proof}

We have proved:

\begin{thm} If \gh\lspace is true in $\mathcal{D}$, then $R$ is von Neumann regular. \label{ghimpliesvnr}\end{thm}

\begin{remark}  In Remark \ref{evendegree}, we noted that if \gh\lspace is true, then the Nilpotence Criterion holds when $\pi_*S$ is concentrated in even degrees.  Since the Annihilator Criterion also holds in this case, we can conclude, as above, that $(x^2) = (x)$ for all $x \in \pi_*S$.  Therefore, \gh\lspace is false in $\mathcal{D}$ if $S$ is connective and $\pi_*S$ is concentrated in even degrees but not concentrated in degree zero. \label{evendegree1}\end{remark}

In order to prove the converse of Theorem \ref{ghimpliesvnr}, we first give circumstances under which $\thicksub{S}$ is exactly the collection of retracts of finite wedges of suspensions of $S$.  In this situation, \gh\lspace is trivially true.  For the next proposition, it is not necessary to assume that $\pi_* S$ is commutative or concentrated in degree zero.

\begin{prop} If every finitely generated submodule of a free $R$-module is a summand, then every element of $\thicksub{S}$ is a retract of a finite wedge of suspensions of $S$.\label{fsub}\end{prop}
\begin{proof} 
Let $\mathcal{C}$ be the collection of objects in $\mathcal{D}$ which satisfy the conclusion of the Proposition.  $\mathcal{C}$ is contained in $\thicksub{S}$ and contains $S$; hence, it suffices to show that $\mathcal{C}$ is thick.  It is clearly closed under suspension and retraction.  To show closure under cofiber sequences, we first consider maps between wedges of spheres.  Let $S_n$ be an $n$-fold wedge of suspensions of $S$ and let $S_m$ be an $m$-fold such wedge.  Consider the cofiber sequence
\begin{equation*}
\xymatrix{\Sigma^{-1} C \ar[r]^-\delta & S_n \ar[r]^-f & S_m \ar[r]^-\phi & C.}
\end{equation*}
By hypothesis, $\im f_*$ is a summand of $\pi_* S_m$, so we get a decomposition $\pi_* S_m \cong \im f_* \oplus \coker f_*$.  Since idempotents split in $\mathcal{D}$ (we assume $\mathcal{D}$ has arbitrary coproducts), we obtain an associated splitting $S_m \simeq K \wdg L$ such that $\pi_* K \cong \im f_*$ and $\pi_* L \cong \coker f_*$.  Since $\im f_*$ is projective, $\ker f_*$ is a summand of $\pi_* S_n$.  Hence we have a decomposition $\pi_* S_n \cong \ker f_* \oplus \coim f_*$ and an associated splitting $S_n \simeq M \wdg N$ with $\pi_*M \cong \coim f_*$ and $\pi_* N \cong \ker f_*$.

Now, $\map{f}{M\wdg N}{K\wdg L}$ has matrix form
\[ \begin{pmatrix} \pi_K f \iota_M & \pi_K f \iota_N \\ \pi_L f \iota_M & \pi_L f \iota_N \end{pmatrix} = \begin{pmatrix} \pi_K f \iota_M & 0 \\ 0 & 0 \end{pmatrix}.\]
All entries but the upper left are zero since they are zero on homotopy groups and have domains that are retracts of wedges of spheres.  We now need to show that $\pi_K f \iota_M$ is an equivalence; if so, then $C \simeq \Sigma N \wdg L$, making $C$ a retract of a wedge of suspensions of the sphere.  Certainly, $\pi_K f \iota_M$ induces an isomorphism of homotopy groups.  This suffices, as our next claim demonstrates.

\vspace{.15in}

\n {\bf Claim.} {\em Let $\map{q}{A}{B}$ be a morphism in $\thicksub{S}$.  If $q$ induces an isomorphism of homotopy groups, then it is an equivalence in $\mathcal{D}$.}
\begin{proof}  First, suppose $M \in \thicksub{S}$ has the property $\pi_* M = [S,M]_* = 0$.  Then the functor $[-,M]_*$ must vanish on all of $\thicksub{S}$.  In particular, $[M,M]_* = 0$, forcing the identity map of $M$ to be trivial.  Hence $M$ is trivial.  Since the cofiber of $q$ satisfies this property, it must be trivial, forcing $q$ to be an equivalence.\end{proof}

Now we consider the general case.  Consider the following commutative diagram, where the rows are cofiber sequences, $W_k$ is a wedge of suspensions of $S$, and the vertical compositions are the identity.
\begin{equation*}
\xymatrix{X \ar[rr]^-f \ar[d]_-{i_1} && Y \ar[rr] \ar[d]^-{i_2} && C \ar[d]^-\psi \\
	W_1 \ar[rr]^-{i_2 f r_1} \ar[d]_-{r_1} && W_2 \ar[rr] \ar[d]^-{r_2} && D \ar[d]^-\phi \\
	X \ar[rr]^-f && Y \ar[rr] && C.}
\end{equation*}
Let $\theta = \phi \psi$.  By the five-lemma, the induced map $\map{\theta^*}{[C,C]_*}{[C,C]_*}$ is an isomorphism; hence, $C$ is a retract of $D$, which is a retract of a wedge of suspensions of $S$ by the above special case.  This completes the proof.
\end{proof}

Next we show that regularity guarantees that the hypothesis of Proposition \ref{fsub} is satisfied. The following proposition follows from \cite[2.7]{goodearl-neumann-1991}.  We include a simple proof for the reader's convenience.

\begin{prop}  If $T$ is a von Neumann regular ring, then every finitely generated submodule of a projective $T$-module is a summand. \label{vnrsummand}\end{prop}
\begin{proof}  It suffices to consider a finitely generated submodule $K$ of $T^n$.  We have a short exact sequence
\[ \xymatrix{0 \ar[r] & K \ar[r] & T^n \ar[r] & Q \ar[r] & 0.} \]
This is a finite presentation of $Q$.  Since $T$ is von Neumann regular, all modules are flat (\cite[4.2.9]{weibel-introduction-1994}).  The module $Q$ is therefore projective by \cite[3.2.7]{weibel-introduction-1994}, so the sequence splits and $K$ is a summand of $T^n$.
\end{proof}

Combining this proposition with Theorem \ref{ghimpliesvnr}, we have proved:

\begin{thm} Let $\mathcal{D}$ be a triangulated category with arbitrary coproducts, and let $S \in \mathcal{D}$ have the property that $R = \pi_*S$ is commutative and concentrated in degree zero.   \gh \lspace is true in $\mathcal{D}$ if and only if $R$ is von Neumann regular.\label{main}\end{thm}

\begin{cor} For commutative rings, \gh\lspace is true in \dr\lspace if and only if $R$ is von Neumann regular. \end{cor}

\begin{remark}  Observe that since Proposition \ref{vnrsummand} applies to noncommuative rings, we may drop the commutativity assumption for one direction of this theorem:  \gh \lspace is true in the derived category of (right) modules over a von Neumann regular ring.\label{noncommincursion}\end{remark}

\begin{remark}  Notice that \gh\lspace is true in $\mathcal{D}$ if and only if every object of $\thicksub{S}$ is a retract of a wedge of suspensions of $S$.  This is a very trivializing set of circumstances.  In \dr, for example, this says that every finite object must be a wedge of suspensions of finitely generated projective modules.  For any connective monogenic stable homotopy category where $\pi_0(S)$ is projective-free, this trivializing condition is equivalent to the requirement that all finite objects be finite wedges of suspensions of the sphere.
\end{remark}

\begin{remark} Continuing the discussion of Remark \ref{evendegree1}, notice that Theorem \ref{main} is true when $\pi_*S$ is concentrated in even degrees; all that must be changed is the definition of regularity, in the obvious way:  a graded ring $R_*$ is regular if for every element $x \in R_k$ there exists an element $y \in R_{-k}$ such that $xyx = x$.  One would also replace `field' with `graded field' in Proposition \ref{vnrchar}.\end{remark}

\begin{example} As a final example, we consider the derived category \dqcoh{X}\lspace of quasi-coherent sheaves of $\mathcal{O}_X$-modules over a scheme $X$, where $\mathcal{O} = \mathcal{O}_X$ is the structure sheaf of $X$.  Since the categories \modcat{R} and \qcoh{X} are equivalent as abelian categories for the affine scheme $X = \spec R$, we have already addressed \gh\lspace in the affine case.  According to \cite[2.6]{hovey-sheaves-2001}, whenever $X$ is a finite-dimensional Noetherian scheme with enough locally frees, \dqcoh{X}\lspace is a unital algebraic stable homotopy category with weak generators the locally free sheaves of finite rank.  The sphere $S$ is the chain complex with $\mathcal{O}$ concentrated in degree zero.  By \cite[III.6.3]{hartshorne-geometry-1977}, $\pi_*S = \Ext_{\qcoh{X}}(\mathcal{O},\mathcal{O}) = H^*(X;\mathcal{O})$.  If $X = \proj{n}_R$, $n$ dimensional projective space over a commutative Noetherian ring $R$, then $\pi_*S = R$ concentrated in degree zero (\cite[III.2.7, III.5.1]{hartshorne-geometry-1977}).  Hence, by Corollary \ref{noethcase}, we conclude that \gh\lspace is true if and only if $R$ is a finite product of fields.  More generally, \gh\lspace is true in \dqcoh{X} for any smooth Fano variety $X$ over a field of characteristic zero; such $X$ have no higher cohomology by the Kodaira Vanishing Theorem (\cite[III.7.15]{hartshorne-geometry-1977}).

We are now able to give an example where \gh\lspace is true for \thicksub{S} but not true if stated for morphisms between any two small objects.  Over $\proj{1}_R$, for example, every locally free sheaf of finite rank is a direct sum of invertible sheaves (\cite[V.2.6]{hartshorne-geometry-1977}), and the invertible sheaves over $\proj{1}_R$ are exactly the twists of the structure sheaf, $\mathcal{O}(k)$ for $k \in \ints$ (\cite[II.6.17]{hartshorne-geometry-1977}).  Hence, the collection of small objects in $\dqcoh{\proj{1}_R}$ is exactly the thick subcategory generated by the $\mathcal{O}(k)$.  Now, $[\mathcal{O},\mathcal{O}(-1)]_* = \Ext_{\qcoh{\proj{1}_R}}(\mathcal{O},\mathcal{O}(-1)) = H^*(\proj{1}_R,\mathcal{O}(-1)) = 0$ (again, use \cite[III.5.1, III.6.3]{hartshorne-geometry-1977}). But the identity map on $\mathcal{O}(-1)$ is nontrivial.\label{dqcoh}\end{example}

\section{Variations}

\n In this section, we change the domain of definition of \gh.  To distinguish among the variants, we will write $\gh_\mathcal{C}$ for the statement:  If $\map{f}{X}{Y}$ is a morphism of objects in $\mathcal{C}$ and $\pi_* f = 0$, then $f \simeq 0$.  Note that $\gh\lspace = \gh_\thicksub{S}$.  We continue to assume that $R = \pi_* S$ is commutative and concentrated in degree zero.

Consider, for example, $\mathcal{C} = \mathcal{D} = \dr$.  Suppose $\gh_\dr \lspace$ is true.  Then \gh\lspace is also true, so $R$ is regular.  Let $A$ and $B$ be $R$-modules with projective resolutions $\chain{P}$ and $\chain{Q}$, respectively.  Now, $\pistar \chain{P} = A$ and $\pistar \chain{Q} = B$ (both concentrated in degree zero), and $[\chain{P}, \chain{Q}]_* = \Ext^*_R(A, B)$.  If $k \geq 0$, then any map $f \in [\chain{P}, \chain{Q}]_k$ must induce the zero map of homotopy groups.  It must therefore be trivial by $\gh_\dr$.  Hence, $\Ext^k_R(A,B) = 0$ for all $k \geq 0$ and all $R$-modules $A$ and $B$.  This forces every $R$-module to be projective, making $R$ Noetherian.  By Proposition \ref{noethcase}, $R$ must be a finite product of fields. Conversely, if $R$ is a finite product of fields, then $\gh_\dr \lspace$ is true (this follows from the fact that Proposition \ref{sums} is also true for $\gh_\dr\lspace$ (the same proof works) and $\gh_\dr \lspace$ is true when $R$ is a field).

For the next variant, we require that $\mathcal{D}$ be a closed symmetric monoidal category (see \cite[\S A.2]{hovey-axiomatic-1997} for all relevant definitions); let $S$ be the unit for the smash product.  We will make use of the product structure and the existence of function objects (and hence a notion of duality).  Corresponding to every finitely generated ideal $I = (x_1, \dots, x_n)$ of $R$, there is a spectrum $S/I = S/x_1 \smsh \cdots \smsh S/x_n$, where $S/x_i$ is the cofiber of $x_i$ as a self map of $S$.  The definition of $S/I$ depends upon the choice of generators for $I$, though $\thicksub{S/I}$ is independent of this choice (\cite[6.0.9]{hovey-axiomatic-1997}).  Motivated by the fact that thick subcategories of this form are involved in the classification of all thick subcategories of $\thicksub{S}$ in \dr, we let $\mathcal{C} = \thicksub{S/I}$ and consider $\gh_\mathcal{C}$.  First, we prove a series of Propositions establishing the relevant homotopy theoretic properties of $S/I$.

Let $I_k = (x_1, \dots, x_k)$.  $S/I_k$ possesses a 'unit' map $\map{\eta_k}{S}{S/I_k}$ defined to be the smash product of the obvious maps $\xymatrix@1{S \ar[r] & S/x_i}$.  Write $\eta$ for $\eta_n$.

\begin{prop} $\pi_0 \eta_k$ is the quotient map $\xymatrix@1{\pi_0 S = R \ar[r] & R/I_k = \pi_0 S/I_k}$, and $\pi_l S/I_k = 0$ for $l < 0$. \label{eta-prop}\end{prop}
\begin{proof} An examination of the cofiber sequence
\[ \xymatrix{S \ar[r]^{x_1} & S \ar[r]^{\eta_1} & S/x_1 \ar[r] & \Sigma S,} \]
coupled with the fact that $\pi_* S$ is concentrated in degree zero, yields the desired result for $k=1$.  Now proceed by induction, using the diagram
\[ \xymatrix{S/I_k \ar[r]^{x_{k+1}} & S/I_k \ar[r] & S/I_{k+1} \ar[r] & \Sigma S/I_k, \\
						   & S \ar[u]^{\eta_k} \ar[ur]_{\eta_{k+1}}} \]
where the triangle commutes and the top row is a cofiber sequence.
\end{proof}

\begin{cor} $z \in I$ if and only if $\eta z = 0$.\label{zinI}\end{cor}

Dually, there are maps $\map{\delta_k}{S/I_k}{\Sigma^k S}$ defined to be the smash product of the obvious maps $\xymatrix@1{S/x_i \ar[r] & \Sigma S}$.  Write $\delta$ for $\delta_n$.  The proof of the following Proposition is dual to the proof of \ref{eta-prop}.

\begin{prop} $\pi_k \delta_k$ is the inclusion $\xymatrix@1{\pi_k S/I_k = \ann I_k \ar[r] & R = \pi_0 S}$, and $\pi_l S/I_k = 0$ for $l > k$.\end{prop}

\begin{cor} $z \in \ann I$ if and only if $z = \delta z'$ for some $z' \in \pi_k S/I$. \end{cor}

\begin{prop} If $\gh_\thicksub{S/I}\lspace$ is true, then $\ann \ann I = I$. \label{doubleannprop}\end{prop}
\begin{proof} Consider $z \in \ann \ann I$.  Since $\im \pi_n \delta = \ann I$, $\pi_* \eta z \delta = 0$.  Hence, $\eta z \delta \simeq 0$ by $\gh_\thicksub{S/I}$.  $\eta$ fits into a cofiber sequence
\[\xymatrix{F \ar[r]^\theta & S \ar[r]^\eta & S/I \ar[r]^\rho & \Sigma F.}\]
Therefore $z\delta = \theta q$ for some $\map{q}{S/I}{F}$.  At this point, we need to make use of the duality functor $D(-) = F(-,S)$, where the objects $F(X,Y)$ are the function objects in the closed symmetric monoidal structure on $\mathcal{D}$.  $S$ is equivalent to its dual, and self maps of $S$ are self dual up to unit.  Further, using the definitions, one can check that $\eta$ and $\delta$ are dual.  Thus $z\delta = \theta q$ implies that $\eta z = \tilde{q} D\theta$, $\map{\tilde{q}}{DF}{S/I}$.  Consider the map $\map{\rho \tilde{q}}{DF}{\Sigma F} \in \pi_0 F\smsh \Sigma F$.  We will prove that $\rho \tilde{q} = 0$ by showing that $\pi_0 F\smsh \Sigma F = 0$.  For each $k$, we have cofiber sequences
\[\xymatrix{F_k \ar[r] & S \ar[r] & S/I_k \ar[r] & \Sigma F_k} \]
which are related via the diagram of cofiber sequences
\[\xymatrix{S \ar[rr]^-{\eta_k} \ar@{=}[d] && S/I_k \ar[d] \ar[r]  & \Sigma F_k \ar[d] \\
		S \ar[rr]^-{\eta_{k+1}} \ar[d] && S/I_{k+1} \ar[r] \ar[d] & \Sigma F_{k+1} \ar[d] \\
		\cdot \ar[rr] && \Sigma S/I_k \ar@{=}[r] & \Sigma S/I_k.} \]
After smashing this diagram with $F$, one can prove, by induction on $k$, that
\[ \pi_l \Sigma F_k \smsh F = \pi_l \Sigma S/I_k = 0 \lspace\text{for}\lspace l \leq 0.\]
In particular, $\pi_0 F \smsh \Sigma F = 0$.  Now, since $\rho \tilde{q} = 0$, we have $\tilde{q} = \eta h$ and $\eta z = \eta h D\theta$.  Since $hD\theta$ and $\theta Dh$ differ by a unit, we have $hD\theta = \theta \tilde{h}$.  Hence $\eta z = \eta \theta \tilde{h} = 0$, and $z \in I$ by Corollary \ref{zinI}.  This completes the proof.
\end{proof}

\begin{cor}  If $\gh_\thicksub{S/I}\lspace$ is true, then $\ann \ann J = J$ whenever $S/J \in \thicksub{S/I}$. \label{doubleann}\end{cor}
\begin{proof}  If $\gh_\thicksub{S/I}\lspace$ is true and $S/J \in \thicksub{S/I}$, then $\gh_\thicksub{S/J}\lspace$ is true.\end{proof}

Our next goal is to prove that this Corollary applies to $I^2$.  The following proposition follows from a straightforward manipulation of cofiber sequences.

\begin{prop} $S/I^2 \in \thicksub{S/I}$. \end{prop}

\n We now prove the main result for this subsection.

\begin{thm} Let $I$ be a finitely generated ideal of $R$ such that $I \cap \ann I = 0$.  Then, $\gh_\thicksub{S/I}$ is true if and only if $R = I \times R/I$ as a direct product of rings and $R/I$ is regular. \end{thm}

\begin{proof}  We prove the `only if' direction first.  Assume $\gh_\thicksub{S/I}\lspace$ is true.  Then, by Corollary \ref{doubleann} and the previous Proposition, $\ann \ann I = I$ and $\ann \ann I^2 = I^2$.  Now take $x \in \ann I^2$ and fix an element $a\in I$.  For all $b \in I$, $xab = 0$ since $ab \in I^2$.  Thus $xa \in I \cap \ann I = 0$; i.e., $x \in \ann I$.  This shows that $\ann I^2 \subseteq \ann I$, which in turn implies that $I \subseteq I^2$.  Therefore $I = I^2$.  By Nakayama's lemma (as stated, for example, in \cite[\S 2.8, Corollary 1]{reid}), there exists an element $r$ in $R$ such that $rI = 0$ and $1-r \in I$.  For $i \in I$, $i = i(1 -r)$; hence, $I$ is a principal idempotent ideal and therefore a ring direct factor of $R$.  The ring direct product decomposition $R \cong I \times T$  corresponds to a splitting $S \simeq A \wdg B$ with $\pi_* A = I$ and $\pi_* B = T$.  Further, $S/I \simeq B \wdg \Sigma B$.  Therefore $\thicksub{S/I} = \thicksub{B}$.  Since there are no maps from $A$ to $B$, $\pi_*(-)$ and $[B, -]_*$ are identical on $\thicksub{B}$.  $T$ is therefore von Neumann regular by Theorem \ref{main}.

Conversely, if $R/I$ is regular and a ring direct factor of $R$, then we obtain, as above, a decomposition $S \simeq A \wdg B$ with $\thicksub{S/I} = \thicksub{B}$ and $[B,B]_* = \pi_* B = R/I$.  By Theorem \ref{main}, $\gh_\thicksub{S/I}\lspace$ is true.
\end{proof}

The next proposition helps us identify situations where the condition $I \cap \ann I = 0$ holds.

\begin{prop}  If $\gh_\thicksub{S/I}\lspace$ is true and $x \smsh S/I \simeq 0$ for all $x \in I$, then $I \cap \ann I = 0$. \end{prop}
\begin{proof}  Let $z \in I \cap \ann I$.  Consider the diagram of cofiber sequences
\[ \xymatrix{F \ar[r]^\theta \ar[d]_0 & S \ar[r]^\eta \ar[d]^z & S/I \ar[r]^\rho \ar@{-->}[d]^h & \Sigma F \ar[d]^0 \\
		  F \ar[r]_\theta & S \ar[r]_\eta & S/I \ar[r]_\rho & \Sigma F,} \] 
where $h$ is any map that makes the diagram commute.  Note that $h$ must factor through both $S$ and $\Sigma F$ ($\eta z \simeq 0$ since $z \in I$).  Since $\pi_0 \Sigma F = 0$, $\pi_* h = 0$.  By $\gh_\thicksub{S/I}$, $h \simeq 0$.  Since $z \in \ann I$, $ z = \delta z'$.  By duality, $z = k \eta$ for some $k \in [S/I, S]_0$.  Since we may take $h = \eta k$ in the above diagram, we have $\eta k \simeq 0$.  Hence $z$ is of the form $\theta l \eta$ for some $l \in [S/I, F]_0$.  We now prove $\theta_*[S/I, F]_0 = 0$, forcing $z$ to be trivial.  By naturality of duality, we have the following commutative diagram, where the vertical maps are duality isomorphisms and the bottom sequence is exact:
\[ \xymatrix{[S/I, F]_0 \ar[rr]^{\theta_*} \ar@{=}[d] && [S/I, S]_0 \ar@{=}[d]&& \\
		  [S, DS/I \smsh F]_0 \ar[rr]_-{(DS/I \smsh \theta)_*} && [S, DS/I \smsh S]_0 \ar[rr]_-{(DS/I \smsh \eta)_*} && [S, DS/I \smsh S/I]_0.} \]
Since $x \smsh S/I \simeq 0$ for all $x \in I$ (and therefore the same is true of the maps $x \smsh DS/I$), the map $DS/I \smsh \eta$ splits.  Hence $\theta_* = 0$, as desired.
\end{proof}

This Proposition applies to the derived category of a ring; it is straightforward to check that $x \smsh S/x \simeq 0$ for all $x \in R$.  We therefore have

\begin{cor}  Let $\mathcal{D} = \dr$, and let $I$ be a finitely generated ideal of $R = \pi_* S$.  $\gh_\thicksub{S/I}\lspace$ is true if and only if $I$ is a ring summand of $R$ with $R/I$ von Neumann regular.\end{cor}

\section{Noetherian Stable Homotopy\label{s-noetherian}}

In this section, we consider the case where  $\pi_* S$ is a Noetherian ring; we do not require that $\pi_* S$ be commutative or concentrated in degree zero.  To begin, we use the following generalization of a result of Freyd's (\cite{me}).  All $\pi_* S$-modules are right modules.

\begin{prop} Let $\mathcal{D}$ be a triangulated category with distinguished object $S$.  Suppose \gh\lspace is true in $\mathcal{D}$.  For every $X \in \thicksub{S}$, if $\pi_* X$ is finitely generated as a graded $\pi_*S$-module, then $X$ is a retract of a finite wedge of suspensions of $S$. \label{fg} \end{prop}

\n We now observe that when $\pi_* S$ is Noetherian, \gh\lspace can only be true under trivial circumstances.

\begin{prop} Let $\mathcal{D}$ be a triangulated category and let  $S$ be an object in $\mathcal{D}$ such that $\pi_*S$ is Noetherian.  Then, \gh\lspace is true in $\mathcal{D}$ if and only if $\thicksub{S}$ is exactly the collection of retracts of finite wedges of suspensions of $S$.\label{noetherian-trivial}\end{prop}
\begin{proof} The `if' direction is trivial.  For the `only if' direction, we show that the collection $\mathcal{C}$ of retracts of finite wedges of suspensions of $S$ is thick.  Since it contains $S$, this gives the desired conclusion.  $\mathcal{C}$ is trivially closed under suspension and retraction.  If $X$ and $Y$ are in $\mathcal{C}$, then $\pi_* X$ and $\pi_* Y$ are finitely generated right modules over a Noetherian ring; hence, the homotopy of any cofiber $C$ of a map from $X$ to $Y$ is also a finitely generated right $\pi_*S$-module.  By Proposition \ref{fg}, this implies that $C$ is a retract of a finite wedge of suspensions of $S$.\end{proof}

\n Further, we record the following related observation concerning the derived category of a ring.

\begin{prop} Let $R$ be a ring, and let \dr\lspace be the derived category of right $R$-modules.  If $\thicksub{S}$ is the collection of retracts of finite wedges of suspensions of $S$, then $R$ is von Neumann regular.\label{dcatobs}\end{prop}
\begin{proof}  Let $x$ be an element of $R$, and let $S/x$ be the cofiber of $\xymatrix@1{S \ar[r]^-x & S}$.  By hypothesis, $\pi_* S/x$ is projective; hence, $\pi_0 S/x = R/(x)_R$ is projective.  Consequently, there is a projection $\xymatrix{R \ar[r] & (x)_R}$ mapping 1 to $xk$ for some $k \in R$, and $xkx = x$.  $R$ is therefore von Neumann regular. \end{proof}

We conjecture that the collection of graded Noetherian rings for which \gh\lspace holds is exactly the collection of semisimple rings.  A ring $R$ is {\em semisimple} if, as a right module over itself, it decomposes as a finite direct sum of simple modules.  It is a fact (\cite[9.4]{anderson-fuller}) that every finitely generated submodule of a free module over a semisimple ring is a summand.  Hence, by Proposition \ref{fsub}, \gh\lspace is true if $\pi_* S$ is semisimple, establishing part of our conjecture.  We will prove the remaining half of the conjecture for the derived category of a ring using the following well-known result (\cite[2.16]{goodearl-neumann-1991}).

\begin{prop} If a ring R is right Noetherian regular, then it is semisimple. \label{noetherian-fpf}\end{prop}
\begin{proof}
A ring is semisimple if and only if every ideal is a direct summand (\cite[9.6]{anderson-fuller}).  Since $R$ is Noetherian, every ideal $I$ of $R$ is finitely generated.  Since $R$ is regular, every finitely generated ideal of $R$ is a summand (\cite[4.2.8]{weibel-introduction-1994}).  Hence,  $R$ is semisimple.
\end{proof}

\begin{cor}  Let $R$ be a Noetherian ring.  Then, \gh\lspace is true in \dr\lspace if and only if $R$ is semisimple.  If $R$ is also commutative, then \gh\lspace is true in \dr\lspace if and only if $R$ is a finite product of fields.\label{noethcase}\end{cor}
\begin{proof} We observed that the 'if' direction is true earlier in this section.  For the 'only if' direction, simply combine Propositions \ref{noetherian-trivial}, \ref{dcatobs}, and \ref{noetherian-fpf}.  The final observation follows from the Artin-Wedderburn theorem, which classifies semisimple rings.\end{proof}

\vspace{.5in}

\bibliographystyle{plain}
\bibliography{gh}

\begin{thebibliography}{10}

\bibitem{anderson-fuller}
Frank~W. Anderson and Kent~R. Fuller.
\newblock {\em Rings and categories of modules}, volume~13 of {\em Graduate
  Texts in Mathematics}.
\newblock Springer-Verlag, New York, second edition, 1992.

\bibitem{bokstedt-neeman}
Marcel B{\"o}kstedt and Amnon Neeman.
\newblock Homotopy limits in triangulated categories.
\newblock {\em Compositio Math.}, 86(2):209--234, 1993.

\bibitem{devinatz-1990}
Ethan~S. Devinatz.
\newblock {$K$}-theory and the generating hypothesis.
\newblock {\em Amer. J. Math.}, 112(5):787--804, 1990.

\bibitem{devinatz-1988}
Ethan~S. Devinatz, Michael~J. Hopkins, and Jeffrey~H. Smith.
\newblock Nilpotence and stable homotopy theory. {I}.
\newblock {\em Ann. of Math. (2)}, 128(2):207--241, 1988.

\bibitem{freyd-1966}
Peter Freyd.
\newblock Stable homotopy.
\newblock In {\em Proc. Conf. Categorical Algebra (La Jolla, Calif., 
1965)},
  pages 121--172. Springer, New York, 1966.

\bibitem{goodearl-neumann-1991}
K.~R. Goodearl.
\newblock {\em von {N}eumann regular rings}.
\newblock Robert E. Krieger Publishing Co. Inc., Malabar, FL, second edition,
  1991.

\bibitem{hartshorne-geometry-1977}
Robin Hartshorne.
\newblock {\em Algebraic geometry}.
\newblock Springer-Verlag, New York, 1977.
\newblock Graduate Texts in Mathematics, No. 52.

\bibitem{hopkins-smith-1998}
Michael~J. Hopkins and Jeffrey~H. Smith.
\newblock Nilpotence and stable homotopy theory. {II}.
\newblock {\em Ann. of Math. (2)}, 148(1):1--49, 1998.

\bibitem{hovey-sheaves-2001}
Mark Hovey.
\newblock Model category structures on chain complexes of sheaves.
\newblock {\em Trans. Amer. Math. Soc.}, 353(6):2441--2457 (electronic), 2001.

\bibitem{hovey-axiomatic-1997}
Mark Hovey, John~H. Palmieri, and Neil~P. Strickland.
\newblock Axiomatic stable homotopy theory.
\newblock {\em Mem. Amer. Math. Soc.}, 128(610):x+114, 1997.

\bibitem{me}
Keir Lockridge.
\newblock The generating hypothesis in general stable homotopy categories.  Ph.D. thesis.  University of Washington, 2006.

\bibitem{reid}
Miles Reid.
\newblock {\em Undergraduate commutative algebra}, volume 29 of {\em London Mathematical Society Student Texts}.
\newblock Cambridge University Press, Cambridge, 1996.

\bibitem{weibel-introduction-1994}
Charles~A. Weibel.
\newblock {\em An introduction to homological algebra}, volume~38 of {\em
  Cambridge Studies in Advanced Mathematics}.
\newblock Cambridge University Press, Cambridge, 1994.

\end{thebibliography}

\end{document}